\def\b#1{{\bf #1}}
\def\i#1{{\it #1}}

\documentstyle[11pt]{article}
\begin{document}

\title{Sequences related to convergents to square root of rationals}
\author{Mario Catalani\\
Department of Economics, University of Torino\\
Via Po 53, 10124 Torino, Italy\\
mario.catalani@unito.it}
\date{}
\maketitle

\bigskip
\bigskip

\section{The Initial Results}
This note has its source in \cite{murthy}, in a different setting.
Consider the system of recurrences, with
$a_0=b_0=1$,
\begin{equation}
\label{eq:prima}
a_n=a_{n-1}+kb_{n-1},
\end{equation}
\begin{equation}
\label{eq:seconda}
b_n=a_{n-1}+b_{n-1}.
\end{equation}
Later on we will generalize and show that the
ratio ${a_n\over b_n}$ is related to convergents to the square root
of rationals.

\noindent
We postpone the proof by
induction of the following formulas ("summation formulas")
\begin{equation}
\label{eq:uno1}
a_{2n}=\sum_{i=0}^n{2n+1\choose 2i+1}k^{n-i},
\end{equation}
\begin{equation}
\label{eq:due2}
a_{2n+1}=\sum_{i=0}^{n+1}{2n+2\choose 2i}k^{n+1-i},
\end{equation}
\begin{equation}
\label{eq:tre3}
b_{2n}=\sum_{i=0}^n{2n+1\choose 2i}k^{n-i},
\end{equation}
\begin{equation}
\label{eq:quattro4}
b_{2n+1}=\sum_{i=0}^n{2n+2\choose 2i+1}k^{n-i}.
\end{equation}
Let $\alpha=k+1+2\sqrt{k},\,\beta=k+1-2\sqrt{k}$. Note that
$\alpha+\beta=2(k+1),\,\alpha\beta=(k-1)^2,\,\alpha-\beta=4\sqrt{k}$.
Using Identities 1.87 and 1.95 in \cite{gould}
and Pascal's Identity (see \cite{koshy})
we have the following
closed form representations
\begin{equation}
\label{eq:unobis}
a_{2n}={\alpha^n+\beta^n\over 2}+\sqrt{k}{\alpha^n-\beta^n\over 2},
\end{equation}
\begin{equation}
\label{eq:duebis}
a_{2n+1}={\alpha^{n+1}+\beta^{n+1}\over 2},
\end{equation}
\begin{equation}
\label{eq:trebis}
b_{2n}={\alpha^n+\beta^n\over 2}+{\alpha^n-\beta^n\over 2\sqrt{k}},
\end{equation}
\begin{equation}
\label{eq:quattrobis}
b_{2n+1}={\alpha^{n+1}-\beta^{n+1}\over 2\sqrt{k}}.
\end{equation}
Let $w(r,\,s)$ denote the recurrence
$$w_n=2(k+1)w_{n-1}-(k-1)^2w_{n-2},\quad w_0=r,\,w_1=s.$$ Furthermore
let $d_n=w(1,\,k+1),\,u_n=w(0,\,2k),\,v_n=w(0,\,2)$. The closed forms are
$$d_n={\alpha^n+\beta^n\over 2},$$
$$u_n={(\alpha-\beta)(\alpha^n-\beta^n)\over 8},$$
$$v_n=2{\alpha^n-\beta^n\over \alpha-\beta}.$$
Note that
\begin{equation}
\label{eq:fondamentale}
kv_n=u_n.
\end{equation}
Then
\begin{equation}
\label{eq:unotris}
a_{2n}=d_n+u_n,
\end{equation}
\begin{equation}
\label{eq:duetris}
a_{2n+1}=d_{n+1},
\end{equation}
\begin{equation}
\label{eq:tretris}
b_{2n}=d_n+v_n,
\end{equation}
\begin{equation}
\label{eq:quattrotris}
b_{2n+1}=v_{n+1}.
\end{equation}
It follows that $a_{2n}=w(1,\,3k+1),\,a_{2n+1}=w(k+1,\,k^2+6k+1),\,
b_{2n}=w(1,\,k+3),\,b_{2n+1}=w(2,\,4k+4)$. Finally
$$a_n=2(k+1)a_{n-2}-(k-1)^2a_{n-4}$$
$$a_0=1,\,a_1=k+1,\,a_2=3k+1,\,
a_3=k^2+6k+1,$$
$$b_n=2(k+1)b_{n-2}-(k-1)^2b_{n-4}$$
$$b_0=1,\,b_1=2,\,b_2=k+3,\,
b_3=4k+4.$$
These fourth order recurrences can be transformed in second order
recurrences in the following way. First of all note that
$$1-2(k+1)x^2+(k-1)^2x^4=\left (1-2x-(k-1)x^2\right )
\left (1+2x-(k-1)x^2\right ).$$
The generating function of $a_n$ is
$${1+(k+1)x+(k-1)x^2-(k-1)^2x^3\over
1-2(k+1)x^2+(k-1)^2x^4}.$$
This can be simplified to
$${1+(k-1)x\over 1-2x-(k-1)x^2},$$
so that we obtain
$$a_n=2a_{n-1}+(k-1)a_{n-2},\quad a_0=1,\,a_1=k+1.$$
Analogously,
the generating function of $b_n$ is
$${1+2x-(k-1)x^2\over
1-2(k+1)x^2+(k-1)^2x^4},$$
and this becomes
$${1\over 1-2x-(k-1)x^2},$$
so that we obtain
$$b_n=2b_{n-1}+(k-1)b_{n-2},\quad b_0=1,\,b_1=2.$$
Writing $\epsilon=1+\sqrt{k},\,\eta=1-\sqrt{k}$ (so that
$\epsilon=\alpha^{1\over 2},\,\eta=\beta^{1\over 2}$) the closed forms are
$$a_n={\epsilon^{n+1}+\eta^{n+1}\over 2},$$
$$b_n={\epsilon^{n+1}-\eta^{n+1}\over 2\sqrt{k}}.$$
Now it is easy to show that
$\lim_{n\longrightarrow\infty}{a_n\over
b_n}=\sqrt{k}$. Indeed
\begin{eqnarray*}
{a_n\over b_n}&=&
{\epsilon^{n+1}+\eta^{n+1}\over 2}\cdot{2\sqrt{k}\over
\epsilon^{n+1}-\eta^{n+1}}\\
&=&\sqrt{k}
{\epsilon^{n+1}+\eta^{n+1}\over
\epsilon^{n+1}-\eta^{n+1}}\\
&=&\sqrt{k}
{1+\left ({\eta\over\epsilon}\right )^{n+1}\over
1-\left ({\eta\over\epsilon}\right )^{n+1}}.
\end{eqnarray*}
Since $\left ({\eta\over\epsilon}\right )^{n+1}$
converges to zero we have that
${a_n\over b_n}$ converges to $\sqrt{k}$.

\medskip

\noindent
Now we are ready to prove by induction the summation formulas.
Assume that Equation~\ref{eq:uno1},
Equation~\ref{eq:due2}, Equation~\ref{eq:tre3} and
Equation~\ref{eq:quattro4} hold
for some integer $n$. This is tantamount to say that
Equation~\ref{eq:unotris},
Equation~\ref{eq:duetris}, Equation~\ref{eq:tretris} and
Equation~\ref{eq:quattrotris} hold
for some integer $n$.
For $n=0$ they are
satisfied since
$$a_0=1={1\choose 1}k^0=d_0+u_0,$$
$$a_1=k+1=\sum_{i=0}^1{2\choose 2i}k^{1-i}=d_1,$$
$$b_0=1={1\choose 1}k^0=d_0+v_0,$$
$$b_1=2={2\choose 1}k^0=v_1.$$
Now
\begin{eqnarray*}
a_{2(n+1)}&=&a_{2n+2}\\
&=&a_{2n+1}+kb_{2n+1}\\
&=&d_{n+1}+kv_{n+1}\\
&=&d_{n+1}+u_{n+1},
\end{eqnarray*}
where the third line is due to the induction hypothesis and we used
Equation~\ref{eq:fondamentale}. Then Equation~\ref{eq:uno1} is
satisfied
through Equation~\ref{eq:unotris}.
\begin{eqnarray*}
b_{2(n+1)}&=&b_{2n+2}\\
&=&a_{2n+1}+b_{2n+1}\\
&=&d_{n+1}+v_{n+1}.
\end{eqnarray*}
Again the third line is due to the induction hypothesis.
Then Equation~\ref{eq:tre3} is
satisfied
through Equation~\ref{eq:tretris}.
\begin{eqnarray*}
a_{2(n+1)+1}&=&a_{2n+3}\\
&=&a_{2n+2}+kb_{2n+2}\\
&=&d_{n+1}+u_{n+1}+kd_{n+1}+kv_{n+1}\\
&=&(1+k)d_{n+1}+2u_{n+1}\\
&=&(1+k){\alpha^{n+1}+\beta^{n+1}\over 2}+{(\alpha^{n+1}-\beta^{n+1}
(\alpha-\beta)\over 4}\\
&=&
{(\alpha+\beta)(\alpha^{n+1}+\beta^{n+1})\over 2}+{(\alpha^{n+1}-\beta^{n+1}
(\alpha-\beta)\over 4}\\
&=&{2\alpha^{n+2}+2\beta^{n+2}\over 4}\\
&=&d_{n+2}.
\end{eqnarray*}
Here the third line is due to the fact that we already proved the
formulas for $a_{2n}$ and $b_{2n}$.
Then Equation~\ref{eq:due2} is
satisfied
through Equation~\ref{eq:duetris}.
For the last relationship consider first of all that
$$u_n+v_n={(\alpha^n-\beta^n)(\alpha+\beta)\over \alpha-\beta}.$$
Then
\begin{eqnarray*}
b_{2(n+1)+1}&=&b_{2n+3}\\
&=&a_{2n+2}+b_{2n+2}\\
&=&2d_{n+1}+u_{n+1}+v_{n+1}\\
&=&\alpha^{n+1}+\beta^{n+1}+
{(\alpha^{n+1}-\beta^{n+1})(\alpha+\beta)\over \alpha-\beta}\\
&=&{(\alpha-\beta)(\alpha^{n+1}+\beta^{n+1})+
(\alpha^{n+1}-\beta^{n+1})(\alpha+\beta)\over \alpha-\beta}\\
&=&2{\alpha^{n+2}-\beta^{n+2}\over \alpha-\beta}\\
&=&v_{n+2}.
\end{eqnarray*}
Here again the third line is due to the formulas already proven.
Then Equation~\ref{eq:quattro4} is
satisfied
through Equation~\ref{eq:quattrotris}.

\section{Other Initial Conditions}
If we start out with $a_0=0,\,b_0=1$ and we denote the resulting sequences
by ${\tilde a_n}$ and ${\tilde b_n}$ we get
$${\tilde a_{2n}}=\sum_{i=0}^{n-1}{2n \choose 2i+1}k^{i+1}=
{\sqrt{k}\over 2}(\alpha^n-\beta^n),$$
$${\tilde a_{2n+1}}=\sum_{i=0}^{n}{2n+1 \choose 2i+1}k^{i+1}=
{\sqrt{k}\over 2}(\alpha^n-\beta^n)+
{k\over 2}(\alpha^n+\beta^n),$$
$${\tilde b_{2n}}=\sum_{i=0}^{n}{2n \choose 2i}k^{i}=
{1\over 2}(\alpha^n+\beta^n),$$
$${\tilde b_{2n+1}}=\sum_{i=0}^{n}{2n+1 \choose 2i}k^{i}=
{\sqrt{k}\over 2}(\alpha^n-\beta^n)+{1\over 2}(\alpha^n+\beta^n).$$
Then
$${\tilde a_n}=2{\tilde a_{n-1}}+(k-1){\tilde a_{n-2}},\quad
{\tilde a_0}=0,\,{\tilde a_1}=k,$$
$${\tilde b_n}=2{\tilde b_{n-1}}+(k-1){\tilde b_{n-2}},\quad
{\tilde b_0}=1,\,{\tilde b_1}=1.$$
It follows easily ($n>0$)
$${\tilde a_n}=kb_{n-1},$$
$${\tilde b_n}=a_{n-1},$$
so that ${{\tilde a_n}\over {\tilde b_n}}\longrightarrow \sqrt{k}$.

\section{Generalization}
The work done allows to analyze easily the following situation. Let us
consider
the sequences
$$u_n=u_{n-1}+kv_{n-1},$$
$$v_n=hu_{n-1}+v_{n-1},$$
with $k$ and $h$ positive integers, where we start with $u_0=1,\,v_0=0$. Then
we get the following summation formulas and closed forms (where
$\alpha_1=1+hk+2\sqrt{hk},\,\beta_1=1+hk-2\sqrt{hk}$)
$$u_{2n}=
\sum_{i=0}^{n}{2n \choose 2i}(hk)^i=
{1\over 2}(\alpha_1^n+\beta_1^n),$$
$$u_{2n+1}=\sum_{i=0}^{n}{2n+1 \choose 2i}(hk)^i=
{\sqrt{hk}\over 2}(\alpha_1^n-\beta_1^n)+
{1\over 2}(\alpha_1^n+\beta_1^n),$$
$$v_{2n}=h\sum_{i=0}^{n-1}{2n \choose 2i+1}(hk)^i=
{\sqrt{h}\over 2\sqrt{k}}(\alpha_1^n-\beta_1^n),$$
$$v_{2n+1}=h\sum_{i=0}^{n}{2n+1 \choose 2i+1}(hk)^i=
{\sqrt{h}\over 2\sqrt{k}}(\alpha_1^n-\beta_1^n)
+{h\over 2}(\alpha_1^n+\beta_1^n).$$
From this, using the same approach as before, we obtain
$$u_n=2u_{n-1}+(hk-1)u_{n-2},\quad
u_0=1,\,u_1=1,$$
$$v_n=2v_{n-1}+(hk-1)v_{n-2},\quad
v_0=0,\,v_1=h.$$
Using the closed forms given before we see that
$${u_{2n}\over v_{2n}}\longrightarrow \sqrt{{k\over h}},\quad
{u_{2n+1}\over v_{2n+1}}\longrightarrow \sqrt{{k\over h}},$$
so that we can conclude
$${u_{n}\over v_{n}}\longrightarrow \sqrt{{k\over h}}.$$

\section{Reduction}
Returning to the initial case,
if $k$ is an odd number the fraction ${a_n\over b_n}$ can be reduced and
new sequences can be defined. More precisely, let $k=2m+1,\,m=0,\,1,\,2,\,
\ldots$. Then
\begin{eqnarray*}
\alpha&=&2m+1+1+2\sqrt{k}\\
&=&2(m+1+\sqrt{k})\\
&=&2\gamma,
\end{eqnarray*}
with $\gamma=m+1+\sqrt{k}$. Analogously, with $\delta=m+1-\sqrt{k}$, we
have
$$\beta=2\delta.$$
Note that $\gamma+\delta=2(m+1),\,\gamma\delta=m^2,\,\gamma-\delta=
2\sqrt{k}=2\sqrt{2m+1}$. Now define
\begin{equation}
c_{2n}=\left ({1\over 2}\right )^na_{2n}={\gamma^n+\delta^n\over 2}
+\sqrt{k}{\gamma^n-\delta^n\over 2},
\end{equation}
\begin{equation}
c_{2n+1}=\left ({1\over 2}\right )^na_{2n+1}=
{\gamma^{n+1}+\delta^{n+1}\over 2},
\end{equation}
\begin{equation}
d_{2n}=\left ({1\over 2}\right )^nb_{2n}={\gamma^n+\delta^n\over 2}+
{\gamma^n-\delta^n\over 2\sqrt{k}},
\end{equation}
\begin{equation}
d_{2n+1}=\left ({1\over 2}\right )^nb_{2n+1}=
{\gamma^{n+1}-\delta^{n+1}\over 2\sqrt{k}}.
\end{equation}
Then of course
$${c_n\over d_n}\longrightarrow \sqrt{k}.$$
Let $u(r,\,s)$ denote the recurrence
$$u_n=2(m+1)u_{n-1}-m^2u_{n-2},\quad u_0=r,\,u_1=s.$$
Then $c_{2n}=u(1,\,3m+2),\,c_{2n+1}=u(m+1,\,m^2+4m+2),\,d_{2n}=u(1,\,m+2),\,
d_{2n+1}=u(1,\,2(m+1))$. And finally
$$c_n=2(m+1)c_{n-2}-m^2c_{n-4},$$
$$d_n=2(m+1)d_{n-2}-m^2d_{n-4},$$
where the initial conditions are determined by the previous recurrences.

\noindent
The generating function of $c_n$ is
$${1+(1+m)x+mx^2-m^2x^3\over 1-2(1+m)x^2+m^2x^4},$$
that of $d_n$ is
$${1+x-mx^2\over 1-2(1+m)x^2+m^2x^4}.$$

\section{Some Identities}
We will work with the sequences
$$a_n=2a_{n-1}+(k-1)a_{n-2},\quad a_0=1,\,a_1=k+1,$$
and
$$b_n=2b_{n-1}+(k-1)b_{n-2},\quad b_0=1,\,b_1=2,$$
with closed forms
$$a_n={\epsilon^{n+1}+\eta^{n+1}\over 2},$$
$$b_n={\epsilon^{n+1}-\eta^{n+1}\over 2\sqrt{k}}.$$
where $\epsilon=1+\sqrt{k},\,\eta=1-\sqrt{k}$.
Using induction and Equations~\ref{eq:prima} and~\ref{eq:seconda}
we get
$$a_n=(k-1)b_{n-1}+b_n,$$
$${a_{n+1}+(k-1)a_{n-1}\over 2k}=b_n,$$
$$kb_n=a_n+(k-1)a_{n-1}.$$
Using the closed forms (and the fact that $\epsilon\eta=1-k$) we get
$$a_n^2-kb_n^2=(1-k)^{n+1},$$
$$a_{2n}=2a_{n-1}a_n-(1-k)^n,$$
$$(k-1)b_{m-1}b_n+b_mb_{n+1}=b_{m+n+1},$$
$$(k-1)a_{m-1}a_n+a_ma_{n+1}=kb_{m+n+1}.$$
Writing $n=m-1$ the last two become
$$(k-1)b_{m-1}^2+b_m^2=b_{2m},$$
$$(k-1)a_{m-1}^2+a_m^2=kb_{2m}.$$
Using these identities we obtain
\begin{equation}
\label{eq:strumento1}
a_{2n}={2ka_n\over k-1}\sqrt{{a_n^2-(1-k)^{n+1}\over k}}
-{2a_n^2\over k-1}-(1-k)^n,
\end{equation}
\begin{equation}
\label{eq:strumento2}
b_{2n}={(1-k)^{n+1}+2kb_n^2-2b_n\sqrt{kb_n^2+(1-k)^{n+1}}\over k-1}.
\end{equation}
Now noting that $a_{2^{n+1}}=a_{2\cdot 2^n}$, writing $2^n$ instead of
$n$ we get recurrences for $a_{2^n}$ and $b_{2^n}$
$$
a_{2^{n+1}}={2ka_{2^n}\over k-1}\sqrt{{a_{2^n}^2-(1-k)^{2^n+1}\over k}}
-{2a_{2^n}^2\over k-1}-(1-k)^{2^n},$$
$$b_{2^{n+1}}={(1-k)^{2^n+1}+2kb_{2^n}^2-2b_{2^n}\sqrt{kb_{2^n}^2+
(1-k)^{2^n+1}}\over k-1}. $$

\section{Newton's Iteration}
The Newton's iteration algorithm (see \cite{wolfram4})
to approximate the square root of integers is given by
$$x_{n+1}={1\over 2}\left (x_n+{k\over x_n}\right ),$$
starting with $x_0=1$. We assume $k>1$. The limit of this recursion is
the fixed point of the mapping
$$f(x)={1\over 2}\left (x+{k\over x}\right ),$$
which is $\sqrt{k}$.

\noindent
For $n\ge 0$ let us write
$$x_n={a_n\over b_n}.$$
We extend the sequences $a_n$ and $b_n$ setting
$a_0=b_0=1$. Then we have
\begin{equation}
\label{eq:uno}
a_n=a_{n-1}^2+kb_{n-1}^2,
\end{equation}
\begin{equation}
\label{eq:due}
b_n=2a_{n-1}b_{n-1}.
\end{equation}
We can generalize considering the recursion
$$x_{n+1}={1\over 2}\left (x_n+{k\over hx_n}\right ),$$
where $h,\,k$ are positive integers with $h\not= k$. This produces
an approximation to $\sqrt{{k\over h}}$. In this case we have
$$a_n=ha_{n-1}^2+kb_{n-1}^2,$$
$$b_n=2ha_{n-1}b_{n-1}.$$

\noindent
We are going to prove using induction that, for $n\ge 2$,
\begin{equation}
\label{eq:tre}
a_n=2a_{n-1}^2-(k-1)^{2^{n-1}}.
\end{equation}
Using Equation~\ref{eq:uno} and Equation~\ref{eq:due} we have $a_1=k+1,\,
b_1=2,\,a_2=(k+1)^2+4k=1+k^2+6k$. If we use Equation~\ref{eq:tre} we have
$a_2=2a_1^2-(k-1)^2=1+k^2+6k$, so that Equation~\ref{eq:tre} is true
for $n=2$. Now assume that it holds for some $n$. Coupled with
Equation~\ref{eq:uno} we get
$$2a_{n-1}^2-w_{n-1}=a_{n-1}^2+kb_{n-1}^2,$$
that is
\begin{equation}
\label{eq:quattro}
a_{n-1}^2=kb_{n-1}^2+w_{n-1},
\end{equation}
where we wrote
$$w_{n-1}=(k-1)^{2^{n-1}}.$$
Now
\begin{eqnarray*}
a_{n+1}&=&a_n^2+kb_n^2\\
&=&(2a_{n-1}^2-w_{n-1})^2+4ka_{n-1}^2b_{n-1}^2\\
&=&4a_{n-1}^4-4a_{n-1}^2w_{n-1}+w_{n-1}^2+4ka_{n-1}^2b_{n-1}^2.
\end{eqnarray*}
On the other hand
\begin{eqnarray*}
2a_n^2-w_n&=&2(a_{n-1}^2+kb_{n-1}^2)^2-w_n\\
&=&2a_{n-1}^4+2k^2b_{n-1}^4+4ka_{n-1}^2b_{n-1}^2-w_n.
\end{eqnarray*}
Now using Equation~\ref{eq:quattro} we have
\begin{eqnarray*}
k^2b_{n-1}^4&=& (a_{n-1}^2-w_{n-1})^2\\
&=&a_{n-1}^4+w_{n-1}^2-2a_{n-1}^2w_{n-1}.
\end{eqnarray*}
Then
\begin{eqnarray*}
2a_n^2-w_n&=&
4a_{n-1}^4-4a_{n-1}^2w_{n-1}+2w_{n-1}^2+4ka_{n-1}^2b_{n-1}^2-w_n\\
&=&a_{n+1},
\end{eqnarray*}
since $w_n=w_{n-1}^2$. This concludes the proof.

\noindent
Using Equation~\ref{eq:quattro} we get
$$a_{n-1}=\sqrt{kb_{n-1}^2+w_{n-1}},$$
so that, through Equation~\ref{eq:due}, we get a recurrence for $b_n$
\begin{equation}
b_n=2b_{n-1}\sqrt{kb_{n-1}^2+w_{n-1}}.
\end{equation}
Now we are going to prove, again by induction, the following closed forms
\begin{equation}
\label{eq:cinque}
a_n={\alpha^{2^n}+\beta^{2^n}\over 2},
\end{equation}
\begin{equation}
\label{eq:sei}
b_n={\alpha^{2^n}-\beta^{2^n}\over 2\sqrt{k}},
\end{equation}
where $\alpha=1+\sqrt{k},\,\beta=1-\sqrt{k}$. Note that $\alpha+\beta=2,\,
\alpha\beta=1-k$.
Hence $a_n$ and $b_n$ are doubly exponential sequences (see \cite{sloane3}).

\noindent
For $n=0$ the closed form for $a_0$ gives ${\alpha+\beta\over 2}=1$; for
$b_0$ gives ${\alpha-\beta\over 2\sqrt{k}}=1$ so Equation~\ref{eq:cinque}
and Equation~\ref{eq:sei} are satisfied. Now assume that
Equation~\ref{eq:cinque} and Equation~\ref{eq:sei} hold for some $n$. Then
\begin{eqnarray*}
a_{n+1}&=&a_n^2+kb_n^2\\
&=&\left ({\alpha^{2^n}+\beta^{2^n}\over 2}\right )^2+
k\left ({\alpha^{2^n}-\beta^{2^n}\over 2\sqrt{k}}\right )^2\\
&=&{\alpha^{2^{n+1}}+\beta^{2^{n+1}}+2(\alpha\beta)^{2^n}\over 4}
+k{\alpha^{2^{n+1}}+\beta^{2^{n+1}}-2(\alpha\beta)^{2^n}\over 4k}\\
&=&{\alpha^{2^{n+1}}+\beta^{2^{n+1}}\over 2}.
\end{eqnarray*}
This concludes the proof for $a_n$.

\noindent
For $b_n$ we have
\begin{eqnarray*}
b_{n+1}&=&2a_nb_n\\
&=&2{\alpha^{2^n}+\beta^{2^n}\over 2}
{\alpha^{2^n}-\beta^{2^n}\over 2\sqrt{k}}\\
&=&{\alpha^{2^{n+1}}-\beta^{2^{n+1}}\over 2\sqrt{k}}.
\end{eqnarray*}
This concludes the proof.

\noindent
Incidentally we have proved, for $n\ge 1$, that $kb_n^2+w_n$ is a perfect
square. With $k=2$, $a_n$ is sequence A001601 and $b_n$ is sequence A051009 in
\cite{sloane}.

\medskip
\noindent
Using Identities 1.87 and 1.95 in \cite{gould}, where $n$ is replaced by
$2^n$, we obtain the following
summation formulas
\begin{equation}
a_n=\sum_{r=0}^{2^{n-1}}{2^n\choose 2r}k^r,
\end{equation}
\begin{equation}
b_n=\sum_{r=0}^{2^{n-1}-1}{2^n\choose 2r+1}k^r,\quad n>0.
\end{equation}
From Equation~\ref{eq:due} we see that
$$b_n=2^n\prod_{r=0}^{n-1}a_r,$$
which implies that $b_n$ is divisible by $2^n$.

\section{Related Sequences}
Let us consider the sequence $c_n$, with $c_0=1,\,c_1=r$, $r>1$, such
that for $n\ge 2$
$$c_n=2c_{n-1}^2-1.$$
From \cite{shallit} we get the closed form
$$c_{n+1}={\gamma^{2^n}+\delta^{2^n}\over 2},$$
where $\gamma=r+\sqrt{r^2-1},\,\delta=r-\sqrt{r^2-1}$. Define the sequence
$d_n$ by
$$d_{n+1}={\gamma^{2^n}-\delta^{2^n}\over 2\sqrt{r^2-1}}.$$
We set $d_0=1$ and we obtain $d_1=1$. Note that $d_2=2r$. Then, for $n\ge 1$,
\begin{eqnarray*}
d_{n+1}&=&2{\gamma^{2^{n-1}}+\delta^{2^{n-1}}\over 2}
{\gamma^{2^{n-1}}-\delta^{2^{n-1}}\over 2\sqrt{r^2-1}}\\
&=&2c_nd_n.
\end{eqnarray*}
Then
$$d_n=2^{n-1}\prod_{i=1}^{n-1}c_i.$$
From this we can evaluate
$$\lim_{n\longrightarrow\infty}\prod_{i=1}^n\left (1+{1\over c_i}\right ).$$
Indeed we have
\begin{eqnarray*}
\prod_{i=1}^n\left (1+{1\over c_i}\right )&=&
\left (1+{1\over c_1}\right )  \left (1+{1\over c_2}\right )\cdots
\left (1+{1\over c_n}\right )\\
&=&{1+c_1\over c_1} {1+c_2\over c_2}\cdots {1+c_n\over c_n}\\
&=&{1+c_1\over c_1} {2c_1^2\over c_2}\cdots {2c_{n-1}^2\over c_n}\\
&=&{(1+c_1)2^{n-1}c_1c_2\cdots c_{n-1}\over c_n}\\
&=&{(r+1)d_n\over c_n}.
\end{eqnarray*}
Using the closed forms of $c_n$ and $d_n$ it is easy to see that
$$\lim_{n\longrightarrow\infty}{c_n\over d_n}=\sqrt{r^2-1}.$$
Hence
$$\lim_{n\longrightarrow\infty}\prod_{i=1}^n\left (1+{1\over c_i}\right )
={r+1\over\sqrt{r^2-1}}=\sqrt{{r+1\over r-1}}.$$
The case $r=3$ is considered in \cite{wilf2}.

\end{document}